\documentclass[a4paper,10pt]{article}
\usepackage[letterpaper,margin=01.1in]{geometry}
\usepackage[utf8]{inputenc}
\usepackage{amssymb}
\usepackage{amssymb, amsmath}
\newtheorem{theorem}{Theorem}
\title{Finite element apriori error estimate for a class of Cardiac Electric models}

\author{Meena Pargaei and B.V. Rathish Kumar\\ Department of Mathematics and Statistics, Indian Institute of Technology, Kanpur
\\ Department of Mathematics and Statistics, Indian Institute of Technology, Kanpur}
\date{}
\begin{document}
\maketitle

\begin{abstract}
In this study we derive the Finite element apriori error estimate for the monodomain cardiac electric model in conjunction with the generic form for a class of nonlinear ionic models.
The analysis establishes a $o(h^2+k)$ space-time convergence. Further the computational realization of the same for different reduced ionic test models is presented.
\end{abstract}

\providecommand{\keywords}[1]
{\textbf{\textit{Keywords:}} #1}

\section{Introduction}
Cardiac electrophysiology is a branch of medicine and biology. Electrophysiological models of heart describe how the electrical potential is generated in every part of the heart. These models consist of coupled ODE-PDE system. 
Electrical activity at the cell level is described by these ODEs and  PDEs describe the same at the tissue level. One of these models is the Bidomain model which is described by two degenerate non-linear parabolic reaction- diffusion
equations which are coupled with a non-linear system of ODE ionic models. The other popular model is the, Monodomain model, which is a simplified form of the Bidomain model. It consists of a non-linear parabolic reaction-diffusion 
equation together with a non-linear system of ODE ionic models. While both finite difference method (FDM)and Finite element method (FEM) have been used to solve these models \cite{fem,fdm1,fdm2,luca}.The space-time
convergence analysis through apriori error estimation has not been reported for these models. 

In this paper we derive the $L^2$ apriori error estimate for the FE analysis of Monodomain model with general form of ionic models. 
In the next section we will describe the cardiac electric models. In section (3) finite element formaulation is presented. In section (4) we derive the $L^2$ apriori
error estimate for semi-discrete and fully discrete system. 
In section (5) a numerical example with space-time convergence results have been presented.

\section{Cardiac Electric Model}
Cardiac tissue is considered as the overlapping of the intra and extracellular continuous domains such that each point in the
intracellular myocardium is also in the extracellular and the two domains are connected with continuous cellular membrane. Bidomain
model(BDM) \cite{luca} is the characterization of such cardiac tissue and it consists of a set of mathematical equations which 
describe the propagation of intra and extracellular electric potentials in cardiac tissue.

The monodomain model(MDM) \cite{luca}, simplified version of the BDM, with the capacity to provide significant information similar
to that of BDM, is used to calculate the action potential. This model consists of a parabolic reaction-diffusion equation coupled
with a system of ODE's which is given as

\begin{equation}
\label{mv}
\frac{\partial v}{\partial t}- div(D(x)\nabla v) - I_{ion}(v,w)=I_{app} \hspace{1.5cm}  \Omega \times (0,T)
\end{equation}
\begin{equation}
\label{mw}
\frac{\partial w}{\partial t}-g(v,w)=0 \hspace{5cm}  \Omega \times (0,T)
\end{equation}

\begin{equation}
\label{in}
v(x,0)= v_0(x,0), \hspace{5mm} w(x,0)=w_0(x,0) \hspace{2cm} \Omega 
\end{equation}

\begin{equation}
\label{bndry}
n^T D(x) \nabla v =0 \hspace{5cm}   \partial \Omega \times (0,T) 
\end{equation}

where $v$ and $w$ are the action potential and gating variables respectively.
$g(v,w)$ and $I_{ion}$ is given by the models at the cell level, called ionic models.

\subsection{Ionic Models}
Hodgkin and Huxely \cite{book} in 1952 gave the first mathematical model that describes accurately the action potential waveform. This model is complex in nature. There are various complex ionic models such as Luo Rudy 1 \cite{luo}, 
Beeler Reuter ,LRd \cite{book}. For large spatial and temporal investigation of any phenomena, various reduced ionic models FitzHugh Nagumo Model(FHNM), Roger-McCulloch Model(RMM) ,Aliev Panfilov Model(APM),
Mitchell Schaeffer Model(MSM) \cite{book,fhn,panfilov,McCulloch,ms} have been described to substantially provide the action potential at low cost.

\section{Galerkin Finite Element Method (GFEM) formulation}
Weak formulation of the system (\ref{mv}-\ref{bndry}) is obtained, find $v(t) \in H^1(\Omega)$ \cite{book}, $w(t) \in L^2(\Omega) or L^{\infty}(\Omega)$\cite{book} and $\psi_1 \in H^1(\Omega) $ , $\psi_2 \in L^2(\Omega)$, 

\begin{align*} 
(v_t,\psi_1) - (\nabla.(D(x)\nabla v), \psi_1)=(I_{ion},\psi_1) \\
	(w_t,\psi_2)=(g,\psi_2)
\end{align*}

Integration by parts together with boundary condition gives
\begin{align} 
\label{weak}
(v_t,\psi_1) - (D(x)\nabla v, \nabla \psi_1)=(I_{ion},\psi_1) \\
(w_t,\psi_2)=(g,\psi_2)
\end{align}

Consider $\{T_h\}_h$  be a member of a family of quasi-uniform triangularizations of $\Omega$ with ${max}_{\tau \in T_h} diam \tau \leq h $ and let $S_h$ be the corresponding finite dimensional space of continuous functions on $\Omega$ which reduces to linear functions in each of the triangles of $T_h$. 

Let $\{\phi_i\}_{i=1}^N$ be the basis functions of $S_h$ so that we can write $v(x,t)= \sum_{i=1}^{N}v_i(t)\phi_i(x)$ and $w(x,t)= \sum_{i=1}^{N}w_i(t)\phi_i(x)$ and equation (\ref{weak}) becomes
\begin{align*}
\Big (\sum_{i=1}^{N}v_{i,t}\phi_i(x),\phi_i \Big)+ \Big (D(x)\nabla \sum_{i=1}^{N}v_i\phi_i,\nabla \phi_i\Big)=(I_{ion},\phi_j) , j=1,2,...,N.
\end{align*}

The matrix form  of the system is given as
\begin{align*}
M v_t+A v= I_{ion}, \hspace{1cm}Mw_t=G,
\end{align*}
where $M_{ij}=((\phi_i,\phi_j))$ , $A_{ij}=(\nabla \phi_i,\nabla \phi_j)$,
$I_{ion}=\Big((I_{ion},\phi_1),...,(I_{ion},\phi_N)\Big)$ and $G=\Big((g,\phi_1),...,(g,\phi_N)\Big)$.

\section{$L^2$ error estimate for the semi-discrete problem}
Let $\Omega$ be a plane convex domain with smooth boundary and consider the problem (\ref{mv}-\ref{bndry}).

Assume that $D(x)$ is symmetric and uniformly positive definite, i.e. there exist $\alpha>0$ such that $\forall x \in R^3 , 
\forall \xi \in R^3$, $\xi^T D(x)\xi\geq \alpha{\mid\xi\mid}^2$.

Introducing the inner product $(\begin{bmatrix} \alpha_1 & \alpha_2 \end{bmatrix}, 
\begin{bmatrix} \beta_1 & \beta_2 \end{bmatrix})_X:=( \alpha_1,\beta_1)+(\alpha_2 , \beta_2)$, and 
the associated norm $\parallel . \parallel_X$, where $(. ,.)$ is the standard $L^2$- inner product.

Let $u={\begin{bmatrix} v & w \end{bmatrix}}^t$ , $ G = \begin{bmatrix} \nabla  & 0 \\ 0 & 0 \end{bmatrix} $ , 
$Gu = G\begin{bmatrix} v \\ w \end{bmatrix} = \begin{bmatrix} \nabla v \\ 0 \end{bmatrix}$.
So,the weak form of the problem becomes
Find $u \in H^1(\Omega) \times L^2(\Omega)$ such that 
\begin{equation}
	\label{wf} 
	(u_t,\psi)_X + (D Gu, G \psi)_X=(F(u),\psi)_X 
\end{equation}
where $F(u)= {\begin{bmatrix} I_{ion}(u) & g(u) \end{bmatrix}}^t$.
Let $\pi_h \colon H^1(\Omega) \times L^2(\Omega) \longrightarrow S_h \times S_h $ be the projection and $u_h={\begin{bmatrix} v_h & w_h \end{bmatrix}}^t$ be the solution of the finite element formulation
\begin{equation}
	\label{fe} 
	(u_{h,t},\psi)_X + (D Gu_h, G \psi)_X=(\pi_h f(u_h),\psi)_X , \forall \psi \in S_h \times S_h 
\end{equation}

\begin{theorem}
	\label{1}
 Let $u$ be the solution of the problem $(\ref{mv}-\ref{in})$, and $u_h$ be the solution of the problem in semi-discrete case. Then, if $I_ion(v,w)$ and $g(v,w)$ is Lipschitz continuous in $v$ and $w$ and $D(x)$ is symmetric and positive definite, we have 
 \begin{align}
 \label{semidiscrete}
 \parallel u_h(T)-u(T) \parallel_X \leq \parallel u_{0,h}-u_0 \parallel_X + C h^2  \parallel u_0 \parallel_X + Ch^2 \int_{0}^{T}(\parallel (F-\pi_h F)(u) \parallel_X+ \parallel u \parallel_X +\parallel u_t \parallel_X )dt
 \end{align}
\end{theorem}
Proof: Decompose the error
\begin{align}
\label{err}
u-u_h = (u - R_hu)+(R_hu-u_h). 
\end{align}

where $R_hu$ is the elliptic projection of the $u$ defined as,\begin{equation}
(D(x) G(R_hu - u) , G\psi)_X =0, \hspace{0.5cm}  \forall \psi \in S_h \times S_h .
\end{equation}
Now we will bound the $\theta= R_hu-u_h$ and $\rho=(u-R_hu)$ separately.
In order to bound $\theta$, note that,
\begin{align*}
(\theta_t,\psi)_X+(D(x)G\theta,G\psi)_X=(F(u)-\pi_h F(u_h),\psi)_X-(\rho_t,\psi)_X,
\end{align*}

Choose $\psi=\theta$,and applying Chauchy-schwartz and lipschitz continuity of F along with boundedness of $L^2$ projection, we get

\begin{align*}
\frac{d}{dt} \parallel \theta \parallel_X \leq \parallel (F-\pi_h F)(u) \parallel_X+ M \parallel u-u_h \parallel_X + \parallel \rho_t \parallel_X.
\end{align*}

After integration we get,
\begin{align}
\label{theta}
\parallel \theta(T) \parallel_X \leq \parallel \theta(0) \parallel_X+ \int_{0}^{T} M(\parallel \theta \parallel_X+\parallel \rho \parallel_X)+\parallel (F-\pi_h F)(u) \parallel_X + \parallel \rho_t \parallel_X dt.
\end{align}

Now apply Gronwall's lemma, and the following bounds for $\parallel \rho \parallel_X$ and $\parallel \rho_t \parallel_X$ taken from the elliptic theory \cite{vidar}
\begin{align}
\label{rho}
\parallel \rho \parallel_X \leq C h^2 \parallel u \parallel_X, 
\parallel \rho_t \parallel_X \leq C h^2 \parallel u_t \parallel_X,
\end{align}
Also, 
\begin{align}
\label{theta0}
\parallel \theta(0) \parallel_X \leq \parallel u_{0,h}-u_0 \parallel_X + \parallel R_hu_0-u_0 \parallel_X \leq \parallel u_{0,h}-u_0 \parallel_X + C h^2  \parallel u_0 \parallel_X
\end{align}
Using (\ref{theta}-\ref{theta0}) in \ref{err}, we arrive at the estimate \ref{semidiscrete}.

\subsection{$L^2$ error estimate for the fully discrete problem}

Let $k$ be the time step, $t_n = nk$, and let $U^n$ be the approximation of $u(t_n)$ in $S_h \times S_h $. We will use backward Euler Galerkin scheme and linearize the problem \ref{wf} by replacing $U^n$ by $U^{n-1}$ to obtain
\begin{equation} 
\label{wft}
(\bar{{\partial}_t} U^n , \psi)_X+ (D(x)GU^n,G\psi)_X = (F(U^{n-1}),\psi)_X, \hspace{1cm} \forall \chi \text{in} S_h.
\end{equation}
where \label{deltU}{$\bar{{\partial}_t} U^n = \frac {1} {k} (U^n-U^{n-1}) $}.

\begin{theorem} 
	\label{2}
	Let $U^n$ and $u$ be solutions of \ref{wft} and \ref{mv}, \ref{mw} respectively. Then, if $I_{ion}(v,w)$ and  $g(v,w)$ are Lipschitz continuous in $v$ and $w$ and $D(x)$ is symmetric and uniformly positive definite, we have
	$\Arrowvert U^n - u(t_n) \Arrowvert_X \leq C \Arrowvert u_{0,h}-u_0 \Arrowvert_X + C(u) (h^2+k) \forall t_n \in \bar{J}$ 
\end{theorem}
\textbf{Proof.} $u^n = u(t_n)$ ,
\begin{equation}
U^n-u^n = (U^n-\widetilde{U}^n) + (U^n-u^n) = \theta ^n + \rho ^n
\end{equation}
where $\widetilde{U}^n$ is the elliptic projection of $u^n$ defined as 

\begin{equation}
(D(x) G(\widetilde{U}^n - u^n) , G\psi)_X = (DG\rho^n , G\psi)_X =0 
\end{equation}

$\rho^n$ will be bounded by Elliptic theory. Now we need to bound only $\theta^n$. For $\psi \in S_h \times S_h$,
\begin{align*}
(\bar{{\partial}_t} \theta^n , \psi)_X+ (D(x)G \theta^n,G\psi)_X = ((\bar{{\partial}_t} U^n , \psi)_X+ (D(x)GU^n,G\psi)_X - (\bar{{\partial}_t} \widetilde{U}^n , \psi)_X
- (D(x)G\widetilde{U}^n,G\psi)_X
\end{align*}
\begin{align*}
=(F(U^{n-1}),\psi)-(u_t^n,\psi)_X-({\overline{{\partial}_t}}\tilde{U}^n-u_t^n,\psi)_X-(D(x)G\tilde{U}^n,G\psi)_X-(D(x)Gu^n,G\psi)_X+(D(x)Gu^n,G\psi)_X
\end{align*}
\begin{align*}
=(F(U^{n-1}),\psi)_X-(F(u^n),\psi)_X-({\overline{{\partial}_t}}(\tilde{U}^n-u^n),\psi)_X
-(D(x)G(\tilde{U}^n-u^n),G\psi)_X-({\overline{{\partial}_t}}u^n-u_t^n,\psi)_X
\end{align*}
Now using Lipschitz continuity of $F$
\begin{align*}
\parallel F(U^{n-1})-F(u^n) \parallel_X \leq C\parallel U^{n-1}-u^n \parallel_X \leq C (\parallel U^{n-1}-u^{n-1} \parallel_X + \parallel u^{n-1}-u^{n} \parallel_X )
\end{align*} 
\begin{equation}
\label{lipf}
\leq C(\parallel \theta^{n-1}  \parallel_X + \parallel \rho \parallel_X +k \parallel {\overline{{\partial}_t}}u^n \parallel_X)
\end{equation}
Take $\psi = \theta^n$ and use the ellipticity of D and (\ref{lipf}) we get,
\begin{align*}
\frac{1}{2} {\overline{{\partial}_t}} {\parallel \theta^{n}  \parallel}_X^2 
+ \mu {\parallel G\theta^{n} \parallel}_X ^2 \leq C (\parallel \theta^{n-1}  \parallel_X + \parallel \rho^{n-1} \parallel_X +k \parallel {\overline{{\partial}_t}}u^n \parallel_X)(\parallel \theta^{n} \parallel_X) 
\end{align*}
\begin{align*}
+ (\parallel {\overline{{\partial}_t}}\rho^n \parallel_X
+\parallel {\overline{{\partial}_t}}u^n-u_t^n \parallel_X)\parallel \theta^{n} \parallel_X
\end{align*}
\begin{align*}
{\overline{{\partial}_t}} {\parallel \theta^{n}  \parallel}_X^2 
\leq C{ (\parallel \theta^{n-1}  \parallel}_X^2 + {\parallel \rho^{n-1} \parallel}_X ^2 +k {\parallel {\overline{{\partial}_t}}u^n \parallel}_X^2) + {\parallel {\overline{{\partial}_t}}\rho^n \parallel_X}^2+ {\parallel {\overline{{\partial}_t}}u^n-u_t^n \parallel_X}^2
+ C {\parallel \theta^{n} \parallel_X}^2 
\end{align*}
\begin{align*}
{\overline{{\partial}_t}} {\parallel \theta^{n}  \parallel}_X^2 
\leq C({\parallel \theta^{n} \parallel}_X^2 + {\parallel \theta^{n-1}\parallel}_X^2+ T_h)
\end{align*}
Using the definition \ref{deltU}
\begin{align*}
{\parallel \theta^{n}  \parallel_X}^2(1-Ck) \leq  {\parallel \theta^{n}  \parallel}_X^2(1+Ck)
+Ck T_h
\end{align*}
\begin{align*}
{\parallel \theta^{n}  \parallel}_X^2 \leq  {\parallel \theta^{n}  \parallel_X}^2(1+Ck)^2
+Ck(1+Ck) T_h
\end{align*}
\begin{align*}
{\parallel \theta^{n}  \parallel}_X^2 \leq (1+Ck)^{n+1} {\parallel \theta^{0}  \parallel}_X^2
+Ck \sum_{j=1}^n (1+Ck)^{n-j+1} T_j
\end{align*}
\begin{align*}
{\parallel \theta^{n}  \parallel}_X^2 \leq C {\parallel \theta^{0}  \parallel}_X^2
+Ck \sum_{j=1}^n T_j
\end{align*}
\begin{align*}
\parallel \rho^j \parallel_X = \parallel{ \tilde U}^j -u^j \parallel_X \leq Ch^2 {\parallel u(jk)\parallel}_X \leq C(u) h^2
\end{align*}
\begin{align*}
\parallel \overline{{\partial}_t} \rho^j \parallel_X = \parallel k^{-1} {\int_{(j-1)k}^{jk}} {\rho_t ds} \parallel_X \leq C(u)h^2
\end{align*}
\begin{align*}
\parallel {\overline{{\partial}_t}}u^n-u_t^n \parallel_X =  \parallel k^{-1} {\int_{(j-1)k}^{jk}} (s-(j-1)k)u_{tt}(s)ds \parallel_X \leq C(u)k
\end{align*}
We have $T_j \leq C(u)(h^2+k)^2$.
Hence,$\parallel \theta^n \parallel_X \leq C \parallel \theta^0 \parallel_X +C(u)(h^2+k)$.

\section{Numerical Test and Discussion}
We consider the monodomain model with the following different ionic models in a square domain $[-1.25,1.25]^2$.
\begin{itemize}
\item FHNM: $I_{ion} = u(u-0.1)(1-u)-w , g(v,w)=u-2w$,
\item RMM: $I_{ion} = u(u-0.1)(1-u)-vw , g(v,w)=u-2w$, 
\item APM: $I_{ion}=-ku(u-a)(u-1)-uw ,$  $g(u,w)=\epsilon' (-ku(u-1-a)-w) $, where $\epsilon'=\epsilon_0 + \mu_1 w/(u+\mu_2)$,

\item MSM: $I_{ion}=-\frac{w}{\tau_{in}}u^2(u-1)-\frac{u}{\tau_{out}},$
				$g(u,w)=  \begin{cases} 
			\frac{1-w}{\tau_{open}} &  u\leq u_{gate},\\
			\frac{-w}{\tau_{close}} & u>u_{gate}.\\
		\end{cases}$
		
	\end{itemize}.

We solved it using linear finite elements in space and Backward Euler in time and the nonlinear terms $I_{ion}$ and $g$ are linearized by taking values at the previous time step. We compute the $L^2$ norm of the error and 
the space and time rate of convergence. 
Initial conditions for $v$ and $w$ are chosen to be 0.2 , 0.1 respectively for all the cases and  $I_{app}=0 ,  dt=h^2 , Dx[n]=h,DT[n]=dt.$
Space rate of convergence(sroc) = $log(L^2(e_{n-1})/L^2(e_n))/log(Dx[n-1])/Dx[n]$, time rate of convergence(troc) = $log(L^2(e_{n-1})/L^2(e_{n})/log(DT[n-1])/DT[n]$, where $e_n : n^{th}$ level error.

\begin{table}[h]
\begin{center}
	\begin{tabular}{ |c|c|c|c|c|c| } 
		\hline
			h & 1/8 & 1/16 & 1/32 &1/64 & 1/128\\ 
		\hline
	FHNM	error & 0.0153718 & 0.00418786 &  0.0010467 & 0.000261422 & 6.53429e-05  \\ 
		\hline
	FHNM	sroc & - & 1.876 & 2.00037 & 2.0014 &2.00027 \\
		\hline
	FHNM	troc & - &0.937999 & 1.00018 & 1.0007 & 1.00014 \\
		\hline
	RMM	error & 0.0293018 &  0.00587156 & 0.00132513 & 0.000321176 &7.93853e-05  \\ 
		\hline
	RMM	sroc & - & 2.31917 & 2.14761 &2.04469 &2.01642  \\
		\hline
	RMM	troc & - &1.15959 & 1.07381 & 1.02235 & 1.00821\\ 
		\hline
	MSM	error & 0.0123688 &  0.0030894 & 0.000772382 & 0.000193005 &4.82713e-05  \\ 
		\hline
        MSM	sroc & - & 2.00131 & 1.99994 &2.0000 &2  \\
		\hline
	MSM	troc & - &1.00065 & 0.999971 & 1.00004 & 0.999998\\
		\hline
	APM	error & 0.0110065 &  0.00273321 & 0.000681779 & 0.000170299 &4.25639e-05  \\ 
		\hline
	APM	sroc & - & 2.00969 & 2.00322 &2.00123 &2.000037  \\
		\hline
	APM	troc & - &1.00484 & 1.00161 & 1.00062 & 1.00018\\
		\hline
	\end{tabular}
	
\caption{$L^2$ norm of the error and the space and time rates of convergence for different ionic models.}
\label{table1}
\end{center}
\end{table} 

From Table \ref{table1} it is clear that as the grid resolution is increased not only the error in $L^2$ norm decreases but also the theortically predicted sroc and troc are realized in all the cases.

\section{Conclusion}
For the Monodomain CEM in a generic framework for ionic models an apriori error estimate under fem approach has been theortically established and computationally verified.

 \section*{Acknowledgement}
 We would like to thank to the DST for the support through Inspire Fellowship.
 
%\section*{References}

\end{document}